\documentclass{amsart}
\usepackage{amssymb,amsthm}
\usepackage{bbm} 
\usepackage{enumerate}
\usepackage{hyperref}
\usepackage[T1]{fontenc} 
\usepackage[latin1]{inputenc} 
\usepackage{ae} 
\usepackage{aecompl} 
\usepackage{mathrsfs}

\newtheorem{theorem}{Theorem}[section]
\newtheorem{lemma}[theorem]{Lemma}
\newtheorem{proposition}[theorem]{Proposition}
\newtheorem{corollary}[theorem]{Corollary}

\newtheorem{conjecture}[theorem]{Conjecture}
\theoremstyle{definition}
\newtheorem{definition}[theorem]{Definition}
\newtheorem{remark}[theorem]{Remark}
\newtheorem{example}[theorem]{Example}

\newcommand{\N}{\mathbb N}
\newcommand{\Z}{\mathbb Z}
\newcommand{\Q}{\mathbb Q}
\newcommand{\R}{\mathbb R}
\newcommand{\C}{\mathbb C}
\newcommand{\K}{\mathbbm k}
\renewcommand{\l}{\ell}
\renewcommand{\i}{\mathbbm i}

\newcommand{\ep}[0]{\varepsilon}
\newcommand{\rh}[0]{\varrho}

\newcommand{\om}[0]{\omega}
\newcommand{\ta}[0]{\tau}
\newcommand{\ph}[0]{\varphi}
\newcommand{\ps}[0]{\psi}

\newcommand{\de}[0]{\delta}
\newcommand{\ze}[0]{\zeta}

\newcommand{\io}[0]{\iota}
\newcommand{\la}[0]{\lambda}
\newcommand{\si}[0]{\sigma}

\newcommand{\La}[0]{\Lambda}

\newcommand{\Ph}[0]{\Phi}
\newcommand{\Ps}[0]{\Psi}

\newcommand{\x}{\bar X}
\newcommand{\ax}{\langle\x\rangle} 
\newcommand{\cx}{[\x]} 
\newcommand{\axy}{\langle X,Y\rangle} 
\newcommand{\cxy}{[X,Y]} 

\renewcommand{\O}{\mathcal O}
\newcommand{\m}{\mathfrak m}
\newcommand{\csim}{\stackrel{\mathrm{cyc}}{\thicksim}}
\DeclareMathOperator{\re}{Re}

\DeclareMathOperator{\tr}{tr}
\DeclareMathOperator{\sym}{Sym}
\DeclareMathOperator{\sign}{sign}

\newcommand{\dd}[2]{\genfrac{}{}{0pt}{}{#1}{#2}}

\title[Connes' embedding conjecture and sums of hermitian squares]
{Connes' embedding conjecture and\\sums of hermitian squares}

\author{Igor Klep}
\address{Igor Klep, Univerza v Ljubljani, Oddelek za matematiko
In\v stituta za matematiko, fiziko in mehaniko,
Jadranska 19, 1111 Ljubljana, Slovénie}
\email{igor.klep@fmf.uni-lj.si}
\thanks{Supported by the Slovenian Research Agency.}
\author{Markus Schweighofer}
\address{Markus Schweighofer, Universit\"at Konstanz,
Fachbereich Mathematik und Statistik,
78457 Konstanz, Allemagne}
\thanks{Supported by the DFG grant ``Barrieren''.}
\email{markus.schweighofer@uni-konstanz.de}
\subjclass[2000]{Primary 11E25, 13J30, 58B34; Secondary 08B20, 47L07, 46L10}
\date{\today}
\keywords{sum of squares, Connes' embedding conjecture, quadratic module,
tracial state, von Neumann algebra}
\begin{document}
\begin{abstract}
We show that Connes' embedding conjecture on von Neumann algebras is
equivalent to the existence of certain algebraic certificates for a
polynomial in noncommuting variables to satisfy the following
nonnegativity condition: The trace is nonnegative whenever
self-adjoint contraction matrices of the same size are substituted
for the variables. These algebraic certificates involve sums of
hermitian squares and commutators. We prove that they always exist
for a similar nonnegativity condition where elements of separable
$\text{II}_1$-factors are considered instead of matrices. Under the
presence of Connes' conjecture, we derive degree bounds for the
certificates.
\end{abstract}
\maketitle

\section{Introduction}

The following has been conjectured in 1976 by Alain Connes
\cite[Section V, pp.~105--107]{con} in his paper on the
classification of injective factors.

\begin{conjecture}[Connes]\label{concon}
If $\om$ is a free ultrafilter on $\N$ and $\mathcal F$ is a
separable $\text{\rm II}_1$-factor, then $\mathcal F$ can be
embedded into the ultrapower $\mathcal R^\om$.
\end{conjecture}

We now explain the notation used in this conjecture. Set
$\N:=\{1,2,3,\dots\}$ and $\N_0:=\{0\}\cup\N$. If $(a_k)_{k\in\N}$
is a sequence in a Hausdorff space $E$ and $\om$ is an ultrafilter
on $\N$, then $\lim_{k\to\om}a_k=a$ means that $\{k\in\N\mid a_k\in
U\}\in\om$ for every neighborhood $U$ of $a$. Such a limit is always
unique and for compact $E$ it always exists. Our reference for von
Neumann algebras is \cite{tak}. When we speak of a trace $\ta$ of a
finite factor $\mathcal F$, we always mean its canonical center
valued trace $\ta:\mathcal F\to\C$ \cite[Definition V.2.7]{tak}.
Such a trace gives rise to the Hilbert-Schmidt norm on $\mathcal F$
given by $\|a\|_2^2:=\ta(a^*a)$ for $a\in\mathcal F$. This norm
induces on $\mathcal F$ a topology which coincides on bounded sets
with the strong operator topology. Let $\mathcal R$ denote the
hyperfinite $\text{\rm II}_1$-factor and $\ta_0$ its trace. Consider
the $C^*$-algebra $\ell^\infty(\mathcal R)
:=\{(a_k)_{k\in\N}\in\mathcal R^\N\mid\sup_{k\in\N}\|a_k\|<\infty\}$
(endowed with the supremum norm). Every ultrafilter $\om$ on $\N$
defines a closed ideal
$I_\om:=\{(a_k)_{k\in\N}\in\ell^\infty(\mathcal
R)\mid\lim_{k\to\om}\|a_k\|_2=0\}$ in $\ell^\infty(\mathcal R)$ and
gives rise to the \textit{ultrapower} $\mathcal
R^\om:=\ell^\infty(\mathcal R)/I_\om$ (the quotient $C^*$-algebra)
which is again a $\text{\rm II}_1$-factor with trace
$\ta_{0,\om}:(a_k)_{k\in\N}+I_\om\mapsto\lim_{k\to\om}\ta_0(a_k)$.
By an \emph{embedding} of $\mathcal F$ into $\mathcal R^\om$, we
always mean a trace preserving $*$-homomorphism.

Recent work of Kirchberg \cite{kir} shows that Connes' conjecture
has several equivalent reformulations in operator algebras and
Banach space theory, among which is the statement that there exists
a unique $C^*$-norm on the tensor product of the universal
$C^*$-algebra of a free group with itself. Voiculescu \cite{voi}
defines a notion of entropy in free probability theory whose
behavior is intimately connected with Connes' conjecture. In this
article, we show that Conjecture \ref{concon} is equivalent to a
purely algebraic statement which resembles recently proved theorems
on sums of squares representations of polynomials. Before presenting
the algebraic reformulation, we need to introduce some notions.

\smallskip
Let always $\K\in\{\R,\C\}$.
As we will rarely need it, we denote the complex imaginary unit by $\i$ so
that the letter $i$ can be used as an index. We denote the complex
conjugate of a complex number $c=a+\i b$ ($a,b\in\R$) by $c^*:=a-\i b$.

We assume that all rings are associative, have a unit element and
that ring homomorphisms preserve the unit element. Throughout the
article, we assume that $n\in\N$ and $\x:=(X_1,\dots,X_n)$ are
variables (or symbols). We write $\ax$ for the monoid freely
generated by $\x$, i.e., $\ax$ consists of \emph{words} in the $n$
letters $X_1,\dots,X_n$ (including the empty word denoted by $1$).
For any commutative ring $R$, let $R\ax$ denote the associative
$R$-algebra freely generated by $\x$, i.e., the elements of $R\ax$
are polynomials in the noncommuting variables $\x$ with coefficients
in $R$. An element of the form $aw$ where $0\neq a\in R$ and
$w\in\ax$ is called a \emph{monomial} and $a$ its
\emph{coefficient}. Hence words are monomials whose coefficient is
$1$. Write $R\ax_k$ for the $R$-submodule consisting of the
polynomials of degree at most $k$ and $\ax_{k}$ for the set of words
$w\in\ax$ of length at most $k$.

\begin{definition}
Let $R$ be a commutative ring. Two polynomials
$f,g\in R\ax$ are called \emph{cyclically equivalent} ($f\csim g$)
if $f-g$ is a sum of commutators in $R\ax$.
\end{definition}

The following remark shows that cyclic equivalence can easily be checked
and that it is ``stable'' under ring extensions in
the following sense: Given an extension of commutative rings $R\subseteq R'$
and $f,g\in R\ax$, then $f\csim g$ in $R\ax$ if and only if $f\csim g$ in
$R'\ax$.

\begin{remark}
Let $R$ be a commutative ring.
\begin{enumerate}[(a)]
\item For $v,w\in\ax$, we have $v\csim w$ if and only if there are
$v_1,v_2\in\ax$ such that $v=v_1v_2$ and $w=v_2v_1$.
\item Two polynomials $f=\sum_{w\in\ax}a_ww$ and $g=\sum_{w\in\ax}b_ww$
($a_w,b_w\in R$) are cyclically equivalent if and only if for each $v\in\ax$,
$$\sum_{\dd{w\in\ax}{w\csim v}}a_w=\sum_{\dd{w\in\ax}{w\csim v}}b_w.$$
\end{enumerate}
\end{remark}

We call a map $a\mapsto a^*$ on a ring $R$ an
\emph{involution} if $(a+b)^*=a^*+b^*$, $(ab)^*=b^*a^*$
and $a^{**}=a$ for all $a,b\in R$. If $*$ is an involution on $R$
(e.g. complex conjugation on $\C$ or the identity on $\R$), then
we extend $*$ to \emph{the} involution on $R\ax$ such that $X_i^*=X_i$.
For each word $w\in\ax$, $w^*$ is its reverse.

\begin{definition}
Let $R$ be a ring with involution $*$. For each subset $S\subseteq R$, we
introduce the set
$$\sym S:=\{g\in S\mid g^*=g\}$$
of its \emph{symmetric elements}. Elements of the form $g^*g$
($g\in R$) are called \emph{hermitian squares}. A subset
$M\subseteq\sym R$ is called a
\emph{quadratic module} if $1\in M$, $M+M\subseteq M$
and $g^*Mg\subseteq M$ for all $g\in R$.
\end{definition}

We can now state the algebraic reformulation of the conjecture.

\begin{conjecture}[Algebraic version of Connes' conjecture]\label{algcon}
Suppose $f\in\K\ax$. If $\K=\R$, assume moreover that $f=f^*$.
Then the following are equivalent:
\begin{enumerate}[{\rm(i)}]
\item $\tr(f(A_1,\dots,A_n))\ge 0$ for all $s\in\N$ and
self-adjoint contractions $A_i\in\K^{s\times s}$;
\label{algcon-matrix}
\item For every $\ep\in\R_{>0}$, $f+\ep$ is cyclically equivalent to an
element in the quadratic module generated by $1-X_i^2$ $(1\le i\le n)$
in $\K\ax$.
\label{algcon-module}
\end{enumerate}
\end{conjecture}

\begin{theorem}\label{mthm}
The following are equivalent:
\begin{enumerate}[\rm (i)]
\item Connes' embedding conjecture {\rm \ref{concon}} holds;
\label{mthm-concon}
\item The algebraic version {\rm \ref{algcon}} of Connes' embedding conjecture
holds;
\label{mthm-algcon}
\item The implication \eqref{algcon-matrix}$\Rightarrow$\eqref{algcon-module}
from Conjecture {\rm \ref{algcon}} {\rm (for $\K=\R$)} holds for all $n\in\N$
and  $f\in\sym\R\ax$.
\label{mthm-matrix}
\end{enumerate}
\end{theorem}

This theorem will be proved in Section \ref{algconsect}. Reformulations of
Connes' conjecture that involve sums of squares have already been given by
Hadwin \cite{had} and  R\u adulescu \cite{r2}. However, Hadwin works with
elements of a certain $C^*$-algebra and R\u adulescu with certain power series
instead of polynomials. In addition, both work with \emph{limits} of sums of
squares. The advantage of our
Conjecture \ref{algcon} is that it is purely algebraic and therefore reveals
the analogy to previously proved theorems on sums of squares representations of
polynomials.

Looking for a counterpart of Conjecture \ref{algcon} for the ring
$\R\cx$ of polynomials in \emph{pairwise commuting} variables, we
replace cyclic equivalence by equality and take the identity
involution. Furthermore, in condition \eqref{algcon-matrix}, the
matrices $A_i$ should now be assumed to commute pairwise. But then
they can be simultaneously diagonalized. One therefore arrives
naturally at the following statement which is a particular case of
Putinar's theorem \cite{put} (we work here over $\K=\R$ since a
complex polynomial which is real on $[-1,1]^n$ has automatically
real coefficients).

\begin{theorem}[Putinar]\label{putinar}
For every $f\in\R\cx$, the following are equivalent:
\begin{enumerate}[{\rm(i)}]
\item $f\ge 0$ on $[-1,1]^n$;
\label{putinar-cube}
\item For all $\ep\in\R_{>0}$, $f+\ep$ lies in the quadratic module
generated by $1-X_i^2$ in $\R\cx$ endowed with the trivial involution.
\label{putinar-module}
\end{enumerate}
\end{theorem}

For noncommuting variables, one can also consider equality instead
of cyclic equivalence. The natural counterpart to Conjecture
\ref{algcon} is then the following particular case of \cite[Theorem
1.2]{hm} (we have omitted the hypothesis $f=f^*$ which is redundant
by \cite[Proposition 2.3]{ks}). For some related results see also
\cite{cim,ks}.

\begin{theorem}[Helton, McCullough]\label{helton}
The following are equivalent for $f\in\K\ax$:
\begin{enumerate}[{\rm(i)}]
\item $f(A_1,\dots,A_n)$ is positive semidefinite for all $s\in\N$ and
self-adjoint contractions $A_i\in\K^{s\times s}$;
\label{helton-matrix}
\item For all $\ep\in\R_{>0}$, $f+\ep$ lies in the quadratic module
generated by $1-X_i^2$ in $\K\ax$. \label{helton-module}
\end{enumerate}
\end{theorem}

The paper is organized as follows. Section \ref{vansect} deals with
polynomials whose trace is not only nonnegative but vanishes. We
prove that these polynomials are sums of commutators. This result is
needed subsequently as a tool. The objective of Section
\ref{algconsect} is to prove Theorem \ref{mthm}. Along the way, we
obtain for example that Conjecture \ref{algcon} holds when matrices
are replaced by elements of $\text{II}_1$-factors (see Theorem
\ref{II_1}). In Section \ref{twovarsect}, we show that Putinar's
Theorem \ref{putinar} implies Conjecture \ref{algcon} for certain
polynomials in two variables. Finally, in Section \ref{bousect} we
establish the existence of certain degree bounds for Conjecture
\ref{algcon}.

\section{Polynomials with vanishing trace}\label{vansect}

\begin{theorem}\label{nullstellensatz}
Let $d\in\N$ and $f\in\K\ax_d$ satisfy
\begin{equation}\label{tr0}
\tr (f(A_1,\dots,A_n))=0
\end{equation}
for all self-adjoint contractions $A_1,\dots,A_n\in\K^{d\times d}$.
In the case $\K=\R$, assume moreover that $f=f^*$. Then $f\csim 0$.
\end{theorem}

\begin{proof}
We call a polynomial $(k_1,\dots,k_n)$-multihomogeneous ($k_i\in\N_0$)
if each of its monomials has for all $i$ degree $k_i$ with respect to the
variable $X_i$. The $(k_1,\dots,k_n)$-multihomogeneous part of a polynomial is
the sum of all its $(k_1,\dots,k_n)$-multi\-homogeneous monomials. Every
polynomial is the sum of its multihomogeneous parts. The
multihomogeneous parts of a symmetric polynomial are symmetric.
We start by proving the following reduction step
which will be used several times during the proof.

{\sc Reduction step.} If $f\in\K\ax$ satisfies \eqref{tr0} for all
self-adjoint contractions $A_1,\dots,A_n\in\K^{d\times d}$, then
all its multihomogeneous parts $g$ satisfy
\begin{equation}\label{tr0mh}
\tr(g(A_1,\dots,A_n))=0
\end{equation}
for all self-adjoint (not necessarily contraction) matrices
$A_1,\dots,A_n\in\K^{d\times d}$.

\emph{Proof of the reduction step.} Fix self-adjoint contractions
$A_1,\dots,A_n\in\K^{d\times d}$. Then for every $\la\in\R$ with
$|\la|\le 1$, the matrix $\la A_1$ is again a self-adjoint
contraction and \eqref{tr0} implies $\tr(f(\la
A_1,A_2,\dots,A_n))=0$. But the latter expression defines a complex
polynomial in $\la$ where the coefficient belonging to $\la^k$ is
$\tr(g_k(A_1,\dots,A_n))$ where $g_k\in\K\ax$ is the sum of all
monomials of $f$ having degree $k$ with respect to $X_1$. Since this
polynomial vanishes at infinitely many points $\la$, all its
coefficients must be zero. This shows that
$\tr(g_k(A_1,\dots,A_n))=0$ for all self-adjoint contractions
$A_1,\dots,A_n\in\K^{d\times d}$. We are therefore reduced to the
case where each $f$ is homogeneous in $X_1$. Now repeat exactly the
same arguments for the other variables. In this way, we see that
\eqref{tr0mh} holds for all multihomogeneous parts $g$ of $f$ and
all self-adjoint contraction matrices $A_i\in\K^{d\times d}$.

\smallskip
As a first application of the now justified reduction step, we see that
our hypothesis implies that \eqref{tr0} holds for all self-adjoint (not
necessarily contraction) matrices. Hence it suffices to show the following
claim for all $k\in\N$ by induction on $k$.

{\sc Claim.} For all $n,d\in\N$ and
$f\in\K\langle X_1,\dots,X_n\rangle_d$ (with $f=f^*$ if $\K=\R$)
having degree at most $k$ in each individual variable $X_i$ and satisfying
\eqref{tr0} for all self-adjoint $A_1,\dots,A_n\in\K^{d\times d}$, we have
$f\csim 0$.

\emph{Induction basis $k=1$.} By the above reduction step and by forgetting
the variables not appearing in $f$, we may assume that $f$ is
$(1,\dots,1)$-homogeneous (also called multilinear), i.e., each variable
appears in each monomial of $f$ exactly once. This means that $f$ can be
written as $f=\sum_{\si\in S_n}a_\si X_{\si(1)}\dotsm X_{\si(n)}$ where
$S_n$ is the symmetric group on $\{1,\dots,n\}$ and $a_\si\in\K$ for all
$\si\in S_n$. By the definition of cyclic equivalence, we have to show that
for each $\ta\in S_n$, the sum over all $a_\si$ such that
$X_{\si(1)}\dotsm X_{\si(n)}$ equals one of the $n$ monomials
$$X_{\ta(1)}\dotsm X_{\ta(n)},\quad
  X_{\ta(2)}\dotsm X_{\ta(n)}X_{\ta(1)},\quad\dots,\quad
X_{\ta(n)}X_{\ta(1)}\dotsm X_{\ta(n-1)}$$
is zero. By renumbering the variables $\x$, we may without loss of generality
assume that $\ta$ is the identity permutation. Let $E_{i,j}\in\K^{d\times d}$
be the matrix with all entries zero except for a one in the $i$-th row and
$j$-th column. Note that $E_{i,j}E_{k,\l}=\de_{j,k}E_{i,\l}$ and
$E_{i,j}+E_{j,i}$ is self-adjoint. Then
it follows from the multilinearity of $f$ that
\begin{align*}
0&=\tr(f(E_{1,2}+E_{2,1},E_{2,3}+E_{3,2},\dots,
         E_{n-1,n}+E_{n,n-1},E_{n,1}+E_{1,n}))\\
&=\tr(f(E_{1,2},E_{2,3},\dots,E_{n-1,n},E_{n,1}))+\dots+
  \tr(f(E_{2,1},E_{3,2},\dots,E_{n,n-1},E_{1,n}))
\end{align*}
where the sum in the last line has $2^n$ terms.
Each of the $2^n-2$ terms represented by the dots must vanish.
This corresponds to the
fact that the only paths on the cyclic graph with $n$ nodes passing through
each of the $n$ edges exactly once are those paths that go through each edge
with the
same orientation (either ``clockwise'' $i\mapsto i+1$ or ``counterclockwise''
$i\mapsto i-1$ modulo $n$). There are only $2n$
such paths which are determined by their starting point and their
orientation. The $n$ clockwise paths show that the first of the $2^n$ terms
is the sum of those $a_\si$ such that $X_{\si(1)}\dotsm X_{\si(n)}$ equals
one of the monomials
\begin{equation}\label{monolist}
X_1\dotsm X_n,\quad
X_2\dotsm X_nX_1,\quad\dots,\quad X_nX_1\dotsm X_{n-1}.
\end{equation}
Calling this sum $a$, we see that $a=0$ is exactly what we have to show.
The $n$ counterclockwise paths show that the last of the $2^n$ terms
is the sum $b$ of those $a_\si$ such that $X_{\si(1)}\dotsm X_{\si(n)}$ equals
one of the monomials
\begin{equation*}
X_n\dotsm X_1,\quad
X_{(n-1)}\dotsm X_1X_n,\quad\dots,\quad X_1X_n\dotsm X_2
\end{equation*}
which are just the monomials arising from \eqref{monolist} by applying the
involution $*$.
Hence $0=a+b$.
In the case $\K=\R$, we use the hypothesis $f=f^*$, to see that $a=b$ and
therefore $a=0$ as desired.
In the case $\K=\C$, additional work is needed.
Choose $\ze\in\C$ such that $\ze^n=\i$. Using similar arguments as above, we
get
\begin{align*}
0&=\tr(f(\ze E_{1,2}+\ze^*E_{2,1},\ze E_{2,3}+\ze^*E_{3,2},\dots,
         \ze E_{n-1,n}+\ze^*E_{n,n-1},\ze E_{n,1}+\ze^*E_{1,n}))\\
&=\ze^n\tr(f(E_{1,2},E_{2,3},\dots,E_{n,1}))+\dots+
  (\ze^*)^n\tr(f(E_{2,1},E_{3,2},\dots,E_{1,n}))\\
&=\i\tr(f(E_{1,2},E_{2,3},\dots,E_{n,1}))-
  \i\tr(f(E_{2,1},E_{3,2},\dots,E_{1,n}))\\
&=\i a-\i b=\i(a-b)
\end{align*}
which together with $a+b=0$ yields $a=0$.

\emph{Induction step from $k-1$ to $k$ {\rm ($k\ge 2$)}.} By the above
reduction step, we can assume that $f$ is $(k_1,\dots,k_n)$-multihomogeneous
where $k_1=\dots=k_m=k$ and $k_i<k$ for all $i\in\{m+1,\dots,n\}$. We assume
$m\ge 1$ since otherwise the induction hypothesis applies immediately.
Now we define recursively a finite sequence $f_0,f_1,\dots,f_m$
of polynomials
$$f_i\in\K\langle X_1,X_1',\dots,X_i,X_i',X_{i+1},X_{i+2}\dots,X_n\rangle$$
by $f_0:=f$ and
\begin{align*}
f_i:=&f_{i-1}(X_1,X_1',X_2,X_2',\dots,X_{i-1},X_{i-1}',X_i+X_i',
              X_{i+1},\dots,X_n)\\
    -&f_{i-1}(X_1,X_1',X_2,X_2',\dots,X_{i-1},X_{i-1}',X_i,
              X_{i+1},\dots,X_n)\\
    -&f_{i-1}(X_1,X_1',X_2,X_2',\dots,X_{i-1},X_{i-1}',X_i',
              X_{i+1},\dots,X_n).
\end{align*}
In other words, each monomial of $f_{i-1}$ gives rise to the $2^k-2$
monomials of $f_i$ which are obtained by replacing at least one but
not all of the occurrences of $X_i$ by $X_i'$. It is important to
note that $f_{i-1}$ can be retrieved from $f_i$ by resubstituting
$X_i'\mapsto X_i$, more exactly
\begin{equation}\label{resubst}
f_{i-1}=\frac 1{2^k-2}
f_i(X_1,X_1',X_2,X_2',\dots,X_{i-1},X_{i-1}',X_i,X_i,X_{i+1},X_{i+2},
    \dots,X_n)
\end{equation}
(we use here that $k\ge 2$). The polynomial $f_m$ has degree at most
$k-1$ with respect to each of its variables and we have
$\tr(f_m(A_1,A_1',\dots,A_m,A_m',A_{m+1},\dots,A_n))=0$ for all
self-adjoint $A_i,A_i'\in\K^{d\times d}$. We now apply the induction
hypothesis (for polynomials in $2m+(n-m)$ variables) to conclude
that $f_m\csim 0$, i.e., $f_m$ is a sum of commutators. Using
\eqref{resubst}, we get successively that $f_{m-1}$, $f_{m-2}$,
$\dots$, $f_0=f$ are also sums of commutators and so $f\csim 0$.
\end{proof}

\begin{remark}\label{indispensable}
For $\K=\R$, the assumption $f=f^*$ in Theorem \ref{nullstellensatz}
is indispensable as shown by $f:=XYZ-ZYX\in\R\langle X,Y,Z\rangle$. For all
$d\in\N$ and all self-adjoint $A,B,C\in\R^{d\times d}$, we have
$\tr(f(A,B,C))=0$ but $f$ is not cyclically equivalent to $0$.
\end{remark}

\begin{proposition}\label{cr-realtrace}
Let $d\in\N$ and $f\in\C\ax_d$ satisfy $\tr(f(A_1,\dots,A_n))\in\R$
for all self-adjoint contractions
$A_1,\dots,A_n\in\C^{d\times d}$. Then there is some $g$ such that
$$f\csim g\in\sym\R\ax_d.$$
\end{proposition}

\begin{proof}
If $f$ were not cyclically equivalent to $p:=\frac{f+f^*}2$, then
$f$ would
not be cyclically equivalent to $f^*$. But then Theorem
\ref{nullstellensatz} would yield complex self-adjoint contraction matrices
$A_i\in\C^{d\times d}$ such that
$$\tr(f(A_1,\dots,A_n))\neq\tr(f^*(A_1,\dots,A_n))=\tr(f(A_1,\dots,A_n))^*,$$
contradicting the hypothesis. Hence $f\csim p$.
Write $p=g+\i h$ with $g,h\in\R\ax$. We have
$g+\i h=p=p^*=(g+\i h)^*=g^*-\i h^*$ and hence $g=g^*$ (and $h=-h^*$).
The ``real trace condition'' which is fulfilled for $f$ by hypothesis, is also
satisfied by $p$ (since $p\csim f$) and $g$ (because
$g=g^*$) and therefore by $\i h$. But this is only possible if
$\tr(h(A_1,\dots,A_n))=0$ for all self-adjoint
$A_i\in\R^{d\times d}$. Applying Theorem \ref{nullstellensatz} again, we obtain
$h\csim 0$. Thus $f\csim g\in\sym\R\ax_d$.
\end{proof}

\section{Algebraic formulation of Connes' conjecture}\label{algconsect}

\begin{definition}\label{tcs}
We call a linear map $\ph:\K\ax\to\K$ a \emph{tracial contraction state} if
\begin{enumerate}[(a)]
\item $\ph(fg)=\ph(gf)$ for all $f,g\in\K\ax$;\label{tcs-t}
\item $|\ph(w)|\le 1$ for all $w\in\ax$;\label{tcs-c}
\item $\ph(f^*f)\ge 0$ for all $f\in\K\ax$;\label{tcs-s}
\item $\ph(1)=1$;\label{tcs-1}
\item (redundant if $\K=\C$, see Remark \ref{redundant} below)
$\ph(f^*)=\ph(f)^*$ for all $f\in\K\ax$.
\label{tcs-*}
\end{enumerate}
\end{definition}

\begin{example}\label{tcs-matrix}
If $A_1,\dots,A_n\in\K^{s\times s}$ are self-adjoint contraction matrices, then
$$\ph:\K\ax\to\K,\qquad f\mapsto\frac 1s\tr(f(A_1,\dots,A_n))$$
is a tracial contraction state.
\end{example}

\begin{remark}\label{redundant}
If $\K=\C$, then \eqref{tcs-*} follows automatically from (a)--(d) in
Definition \ref{tcs}. Indeed, it follows from \eqref{tcs-s} and the identity
\begin{equation}\label{non-useful}
f=\left(\frac{f+1}2\right)^2-\left(\frac{f-1}2\right)^2
\end{equation}
that $\ph(f)\in\R$ for $f\in\sym\C\ax$. Now use that
$\C\ax=\sym\C\ax\oplus\i\sym\C\ax$ as a real vector space which follows from
the identity
\begin{equation}\label{useful}
f=\frac{f+f^*}2+\i\frac{f-f^*}{2\i}.
\end{equation}
\end{remark}

\begin{remark}\label{tcs-eq}
In Definition \ref{tcs}, condition \eqref{tcs-c} can equivalently be replaced
by each of the following conditions:
\begin{enumerate}[(a''')]
\item[(b')] $\ph$ is a contraction with respect to the $1$-norm on $\K\ax$
defined by
$$\Big\|\sum_{w\in\ax}a_ww\Big\|_1:=
\sum_{w\in\ax}|a_w|\qquad(\text{$a_w\in\K$, only finitely many $\neq 0$});$$
\label{tcs-eq-1}
\item[(b'')] The set $\{\ph(X_i^{2k})\mid k\in\N,1\le i\le n\}$ is bounded;
\label{tcs-eq-2k}
\item[(b''')] $\lim\inf_{k\to\infty}|\ph(X_i^{2k})|<\infty$ for
$i\in\{1,\dots,n\}$.
\label{tcs-eq-lim}
\end{enumerate}
For details, consult \cite[Theorem 1.3]{had}.
\end{remark}

\begin{definition}
For any commutative ring $R$ with involution, we denote by
$M_R^{(n)}\subseteq\sym R\ax$ the quadratic module generated by
$1-X_1^2,\dots,1-X_n^2$ in $R\ax$. Most of the time, there will be
no doubt about the number $n$ of variables and we will simply write
$M_R$ instead of $M_R^{(n)}$.
\end{definition}

\begin{remark}\label{identity}
In any $\Q$-algebra $R$, the identity
$$1-a+\frac 1ma^m=\frac 1m+\frac 1m(1-a)^2\sum_{k=0}^{m-2}(m-1-k)a^k$$
holds for all $m\in\N$ and $a\in R$.
\end{remark}

\begin{lemma}\label{tcs-module}
In Definition {\rm \ref{tcs}}, conditions \eqref{tcs-c} and \eqref{tcs-s} can
be replaced by the condition $\ph(M_\K)\subseteq\R_{\ge 0}$.
\end{lemma}

\begin{proof}
Assume that $\ph(M_\K)\subseteq\R_{\ge 0}$. Condition \eqref{tcs-s}
follows immediately since the set of all hermitian squares is
contained in $M_\K$. For $w\in\ax$, $\mu\in\K$ with $|\mu|=1$,
$s\in\N$ and self-adjoint contraction matrices
$A_1,\dots,A_n\in\K^{s\times s}$,
$$\Big(1-\frac{\mu w+(\mu w)^*}2\Big)(A_1,\dots,A_n)$$
is positive semidefinite. Hence by Theorem \ref{helton},
$1-\frac{\mu w+(\mu w)^*}2+\ep\in M_\K$ for every
$\ep\in\R_{>0}$. This implies $\ph(1-\frac{\mu w+(\mu w)^*}2)\ge 0$ and
so $\re(\mu\ph(w))=\re\ph(\mu w)\le 1$. Since $\mu\in\K$ with $|\mu|=1$
was arbitrary, this implies $|\ph(w)|\le 1$.

For the converse, let $g\in\K\ax$ be arbitrary. Then for every
$m\in\N$,
\begin{align*}
g^*(1-X_i^2)g&= g^*(1-X_i^2+\frac 1mX_i^{2m})g - \frac 1mg^*X_i^{2m}g \\
&= g^*
\Big(\frac 1m + \frac 1m (1-X_i^2)^2\sum_{k=0}^{m-2}(m-1-k)X_i^{2k}\Big)g-
\frac 1mg^*X_i^{2m}g
\end{align*}
by Remark \ref{identity}.
By applying $\ph$ to the last expression, the first summand becomes
nonnegative by
\eqref{tcs-s}, while $\frac 1m\ph(g^*X_i^{2m}g)$ goes to zero when $m\to\infty$
since $\ph$ is continuous with respect to the $1$-norm by \eqref{tcs-c}.
This proves that
$\ph(g^*(1-X_i^2)g)\ge 0$. Hence $\ph(M_\K)\subseteq\R_{\ge 0}$.
\end{proof}

\begin{definition} If $R$ is a ring with involution $*$ and
$M\subseteq\sym R$ is a quadratic module, then we define its
\emph{ring of bounded elements}
$$H(M):=\{g\in R\mid N-g^*g \in M\ \text{for some $N\in\N$}\}.$$
This is indeed a $*$-subring of $R$ as proved in \cite[Lemma 4]{vid}.
\end{definition}

In algebra, one says that a quadratic module $M\subseteq\sym R$ is
archimedean if $H(M)=R$. Unfortunately, this has a completely
different meaning in the context of ordered vector spaces
\cite[p.~202, \S 22A]{hol}. We avoid this terminology and instead
use the concept of algebraic interior (or core) points \cite[p.~7,
\S 2C]{hol}.

\begin{definition}
Let $V$ be a $\K$-vector space and $C\subseteq V$.
A vector $v\in V$ is called an
\emph{algebraic interior point} of $C$ if for each $u\in V$ there is some
$\ep\in\R_{>0}$ such that $v+\la u\in C$ for all $\la\in\R$ with
$0\le\la\le\ep$.
\end{definition}

The following is well-known but so important for us that we give a proof of it.

\begin{proposition}\label{algint}
If $R$ is an $\R$-algebra and $M\subseteq\sym R$ a quadratic module, then
$H(M)=R$ if and only if $1$ is an algebraic interior point of $M$ in $\sym R$.
\end{proposition}

\begin{proof}
If $1$ is an algebraic interior point of $M$ in $\sym R$ and $g\in R$, we find
some $N\in\N$ such that $1-\frac 1N g^*g\in M$, i.e., $N-g^*g\in M$.

Conversely, suppose that $H(M)=R$ and let $u\in\sym R$ be given. Then
$u=\left(\frac{u+1}2\right)^2-\left(\frac{u-1}2\right)^2$. Choose $N\in\N$
such that $N-\left(\frac{u-1}2\right)^2\in M$ and set $\ep:=\frac 1N$.
Then $1+\la u\in M$ for all $\la\in\R$ with $0\le\la\le\ep$.
\end{proof}

\begin{lemma}\label{algint1}
If $R$ is a $*$-subfield of $\C$, then $H(M_R)=R\ax$.
\end{lemma}

\begin{proof}
We have $R\subseteq H(M_R)$ and $1-X_i^2\in M_R$, hence $X_i\in H(M_R)$.
Since $H(M_R)$ is a subring of $R\ax$, this implies $H(M_R)=R\ax$.
\end{proof}

\begin{theorem}\label{II_1}
For $f\in\C\ax$, the following are equivalent:
\begin{enumerate}[{\rm (i)}]
\item $\ta(f(A_1,\dots,A_n))\ge 0$ for every separable $\text{\rm II}_1$-factor
$\mathcal F$
with trace $\ta$ and all self-adjoint contractions
$A_1,\dots,A_n\in \mathcal F$;
\label{II_1-trace}
\item $\ph(f)\ge 0$ for all tracial contraction states $\ph$ on $\C\ax$;
\label{II_1-state}
\item For every $\ep\in\R_{>0}$, $f+\ep$ is cyclically equivalent to an element
of $M_\C$.
\label{II_1-module}
\end{enumerate}
\end{theorem}

\begin{proof}
It is immediate from Lemma \ref{tcs-module} that \eqref{II_1-module} implies
\eqref{II_1-state}. It is trivial that \eqref{II_1-state} implies
\eqref{II_1-trace}. To see that \eqref{II_1-trace} implies
\eqref{II_1-module}, we proceed as follows.
Suppose that there is $\ep>0$ such that $f+\ep$ is
not cyclically equivalent to an element of $M_\C$. We start by
constructing a tracial contraction state $L$ on $\C\ax$ such that
$L(f)\notin\R$ or $L(f)<0$.

If $f$ is not cyclically equivalent to
\emph{any} symmetric element, then Proposition \ref{cr-realtrace}
yields a tracial contraction state $L:\C\ax\to\C$ coming from
matrices (cf. Example \ref{tcs-matrix}) such that $L(f)\notin\R$.

If $f$ is cyclically equivalent to a symmetric element of $\C\ax$,
then we may assume without loss of generality that $f$ \emph{is}
symmetric. Define $U:=\{g\in\sym\C\ax\mid g\csim 0\}$. Then $M_\C+U$
is a convex cone in the real vector space $\sym\C\ax$. By Lemma
\ref{algint1}, $1$ is an algebraic interior point of $M_\C$ and
therefore of $M_\C+U$. Since $f+\ep\notin M_\C+U$ and $M_\C+U$
possesses an algebraic interior point, we can apply the
Eidelheit-Kakutani separation theorem \cite[p.~15, \S 4B
Corollary]{hol} to obtain an $\R$-linear functional
$L_0:\sym\C\ax\to\R$ such that $L_0(M_\C+U)\subseteq\R_{\ge 0}$ and
$L_0(f+\ep)\in \R_{\le 0}$. In particular, $L_0(U)=\{0\}$. Using
\eqref{useful}, $L_0$ can be extended uniquely to a $\C$-linear
functional $L$ on $\C\ax$. Obviously, $L$ is a state. To prove that
$L$ is tracial, let $g,h\in\C\ax$ be arbitrary and write $g=g_1+\i
g_2$ and $h=h_1+\i h_2$ for $g_1,g_2,h_1,h_2\in\sym\C\ax$. Then
$[g,h]=[g_1,h_1]+\i [g_2,h_1]+\i [g_1,h_2]- [g_2,h_2]$. The second
and the third summand are symmetric commutators and are thus mapped
to $0$ by $L$. Similarly, $L([g_j,h_j])=-\i L([\i g_j,h_j]))=0$ for
$j=1,2$. Thus $L([g,h])=0$, as desired.

In both cases we obtain a tracial contraction state $L$ with
$L(f)\notin\R_{\ge 0}$. (Note that this already proves
$\eqref{II_1-state} \Rightarrow \eqref{II_1-module}$.)

Endow $\C\ax$ with the $1$-norm defined in Remark \ref{tcs-eq}. By
the Banach-Alaoglu theorem \cite[p.~70, \S 12D Corollary 1]{hol},
the convex set of all tracial contraction states is weak
$*$-compact. Thus by the Krein-Milman theorem \cite[p.~74, \S 13B
Theorem]{hol} we may assume that $L$ is an extreme tracial
contraction state.

We now apply the Gelfand-Naimark-Segal construction with $L$.
By the Cauchy-Schwarz inequality for semi-scalar products,
$N:=\{p \in \C\ax \mid L(p^*p)=0\}$ is a subspace of $\C\ax$.
Similarly, we see that
\begin{equation}\label{gns}
\langle\overline p,\overline q\rangle:=L(q^\ast p)
\end{equation}
defines a scalar product on $\C\ax/N$, where
$\overline p:=p+N$ denotes the residue class of $p\in\C\ax$ modulo $N$.
Let $E$ denote the completion of $\C\ax/N$ with respect to this scalar product.
Since $1\notin N$, $E$ is nontrivial. Observe that $E$ is separable.

To prove that $N$ is a left ideal of $\C\ax$, we fix $i\in\{1,\dots,n\}$ and
show that $X_iN\subseteq N$.
Since $1-X_i^2\in M_\C$ for every $i$, we have
\begin{equation}\label{key-dot-product-relation}
0 \leq L(p^*X_i^2p) \leq L(p^* p)
\end{equation}
for all $p \in \C\ax$. Hence $L(p^*X_i^2p)=0$ for all
$p \in N$, i.e., $X_ip \in N$.

Because $N$ is a left ideal, the map
$$\La_i:\C\ax/N\to\C\ax/N,\;\overline p\mapsto \overline{X_ip}$$
is well-defined for each $i$. Obviously, it is linear and it is self-adjoint
by the definition (\ref{gns}) of the scalar product.
By $(\ref{key-dot-product-relation})$,
$\La_i$ is bounded with norm $\leq 1$ and thus extends to a self-adjoint
contraction $\hat X_i$ on $E$.

Let $\mathcal F$ denote the von Neumann subalgebra of $\mathcal
B(E)$ generated by $\hat X_1,\dots,\hat X_n$ and let $\ta$ denote
the mapping
\begin{equation}\label{trace}
\sum_w a_w \hat w
\mapsto \big\langle \sum_w a_w \hat w(1),1\big\rangle=
L\big(\sum_w a_w w\big).
\end{equation}
$\ta$ is easily seen to be a tracial state on the algebra generated
by $\hat X_1,\dots,\hat X_n$. By continuity, $\ta$ extends uniquely
to a faithful tracial state on $\mathcal F$. Moreover, $1$ is a
separating vector for $\ta$. Hence $\mathcal F$ is a finite von
Neumann algebra \cite[Theorem V.2.4]{tak} and thus can be decomposed
as $\mathcal F=\mathcal F_{\rm I} \oplus \mathcal F_{\rm II}$, where
$\mathcal F_{\rm I}$ and $\mathcal F_{\rm II}$ are finite von
Neumann algebras of type I, respectively II \cite[Theorem
V.1.19]{tak}. Since $L$ was an extremal tracial contraction state,
we have $\mathcal F_{\rm I}=\{0\}$ or $\mathcal F_{\rm II}=\{0\}$.
Assume that the latter holds. Then $\mathcal F$ is a finite type I
von Neumann algebra, hence of type I$_n$ for some $n\in\N$ and is
isomorphic to $n\times n$ matrices over its center \cite[Theorem
V.1.27]{tak}. By \eqref{trace}, $1$ is a trace vector for $\ta$, so
$n=1$, i.e., $\mathcal F$ is abelian. Since $E$ is separable,
$\mathcal F$ can be written as a direct integral of I$_1$-factors
(i.e., $\C$) \cite[Theorem IV.8.21]{tak}. From this decomposition it
follows by assumption \eqref{II_1-trace} that $\ta(\hat f)\geq 0$.
But $\ta(\hat f)=L(f)\notin\R_{\geq 0}$, contradiction.

Hence we may assume that $\mathcal F$ is a type II$_1$ von Neumann
algebra with trace $\ta$. As above, write $\mathcal F$ as a direct
integral of II$_1$-factors and $\ta$ as a direct integral of
(faithful) tracial states. It follows from assumption
\eqref{II_1-trace} that $\ta(\hat f)\geq 0$, again a contradiction
to $\ta(\hat f)=L(f)\notin\R_{\geq 0}$.
\end{proof}

\begin{lemma}\label{cr-intersection}
$M_\C\cap\R\ax=M_\R$. Moreover, if $f\in\R\ax$ is cyclically equivalent to an
element of $M_\C$, then it is cyclically equivalent to an
element of $M_\R$.
\end{lemma}

\begin{proof}
Set $g_0:=1$ and $g_i:=1-X_i^2$ for $i\in\{1,\dots,n\}$ and suppose that
$$\sum_{i=0}^n\sum_j(p_{ij}+\i q_{ij})^*g_i(p_{ij}+\i q_{ij})\in\R\ax$$
where $p_{ij},q_{ij}\in\R\ax$. We have to show that this sum lies in $M_\R$.
Since it lies in $\R\ax$, it is enough to show that it lies in $M_\R$ after
adding its complex conjugate (which is the sum itself). But this is even true
for each particular term in the sum since
$$(p_{ij}+\i q_{ij})^*g_i(p_{ij}+\i q_{ij})+
  (p_{ij}-\i q_{ij})^*g_i(p_{ij}-\i q_{ij})=
  2(p_{ij}^*g_ip_{ij}+q_{ij}^*g_iq_{ij})\in M_\R.$$
For the second statement, let $f+\sum_{i=1}^t [g_{i1},g_{i2}]+\i\sum_{i=1}^t
[h_{i1},h_{i2}]\in M_\C$ for $g_{ij},h_{ij}\in\R\ax$. By applying the complex
conjugation and adding both equations, we obtain
$f+\sum_{i=1}^t [g_{i1},g_{i2}]\in M_\C\cap\R\ax=M_\R$.
\end{proof}

The polynomial from Remark \ref{indispensable} shows that the assumption
$f=f^*$ cannot be omitted in the next two lemmas.

\begin{lemma}\label{rc}
For $f\in\sym\R\ax$, the following are equivalent:
\begin{enumerate}[{\rm (i)}]
\item $\ph(f)\ge 0$ for all tracial contraction states $\ph$ on $\R\ax$;
\label{rc-stater}
\item $\ph(f)\ge 0$ for all tracial contraction states $\ph$ on $\C\ax$.
\label{rc-statec}
\end{enumerate}
\end{lemma}

\begin{proof}
If \eqref{rc-statec} holds and $\ep\in\R_{>0}$, then $f+\ep$ is
cyclically equivalent to an element of $M_\C$ by the implication
\eqref{II_1-state} $\Rightarrow$ \eqref{II_1-module} in Theorem
\ref{II_1}. Hence it is cyclically equivalent to an element of
$M_\R$ by Lemma \ref{cr-intersection} and so $\ph(f)\geq 0$ for all
tracial contraction states $\ph$ on $\R\ax$ by Lemma
\ref{tcs-module}. Conversely, suppose that \eqref{rc-stater} holds
and let $\ph$ be a tracial contraction state on $\C\ax$. Then
$$\ps:\R\ax\to\R,\ p\mapsto\frac{\ph(p)+\ph(p)^*}2$$
is a tracial contraction state. Therefore $\ph(f)=\ps(f)\geq 0$.
\end{proof}

\begin{lemma}\label{cr}
For $f\in\sym\R\ax$, the following are equivalent:
\begin{enumerate}[\rm (i)]
\item $\tr(f(A_1,\dots,A_n))\ge 0$ for all $s\in\N$ and
self-adjoint $A_i\in\R^{s\times s}$;
\label{cr-r}
\item $\tr(f(A_1,\dots,A_n))\ge 0$ for all $s\in\N$ and
self-adjoint $A_i\in\C^{s\times s}$.
\label{cr-c}
\end{enumerate}
\end{lemma}

\begin{proof}
It is trivial that \eqref{cr-c} implies \eqref{cr-r}. For the other
implication, we use the usual identification of a complex number $a+\i b$
($a,b\in\R$) with the real matrix $$\begin{pmatrix}a&-b\\b&a\end{pmatrix}.$$
Every self-adjoint complex matrix defines in this way
a self-adjoint real matrix of double size with double trace.
We leave the details to the reader.
\end{proof}

\begin{corollary}\label{contr}
For $f\in\sym\R\ax$, the following are equivalent:
\begin{enumerate}[{\rm (i)}]
\item $\ph(f)\ge 0$ for all tracial contraction states $\ph$ on $\R\ax$;
\label{contr-state}
\item For every $\ep\in\R_{>0}$, $f+\ep$ is cyclically equivalent to an element
of $M_\R$.
\label{contr-module}
\end{enumerate}
\end{corollary}

\begin{proof}
The implication \eqref{contr-state} $\Rightarrow$ \eqref{contr-module}
follows from Lemma \ref{rc}, Theorem \ref{II_1} and Lemma
\ref{cr-intersection}, while the converse follows from Lemma \ref{tcs-module}.
\end{proof}

The equivalence of \eqref{allultra}, \eqref{oneultra} and
\eqref{approximate} in the next theorem is well-known
\cite{had,r1,r2}. With condition \eqref{approximate}, one can
reformulate Connes' Conjecture \ref{concon}
without recourse to ultraproducts. Our contribution is the new
condition \eqref{nonplusultra}. The implications
\eqref{allultra}$\Rightarrow$\eqref{oneultra}$\Rightarrow$\eqref{nonplusultra}
are easy. The proof of
\eqref{nonplusultra}$\Rightarrow$\eqref{approximate} uses arguments
similar to those of Hadwin \cite[p.~1789]{had} and R\u adulescu
\cite[p.~232]{r1}. Since we work with polynomials, we can even argue
in a simpler way and therefore include a proof. For the sake of
completeness, we also include an elementary proof of
\eqref{approximate}$\Rightarrow$\eqref{allultra} which resembles the
proof of \cite[Lemma 5.22]{con}.

\begin{proposition}\label{approx}
For every separable $\text{\rm II}_1$-factor $\mathcal F$ with trace
$\ta$, the following are equivalent:
\begin{enumerate}[\rm(i)]
\item\label{allultra}
For every free ultrafilter $\om$ on $\N$, $\mathcal F$ is embeddable
in $\mathcal R^\omega$;
\item\label{oneultra}
There is an ultrafilter $\om$ on $\N$ such that $\mathcal F$ is
embeddable in $\mathcal R^\omega$;
\item\label{nonplusultra}
For each $n\in\N$ and $f\in\C\ax$, condition \eqref{algcon-matrix}
from Conjecture {\rm\ref{algcon}} implies $\ta(f(A_1,\dots,A_n))\ge
0$ for all self-adjoint contractions $A_1,\dots,A_n\in\mathcal F$;
\item\label{approximate}
For all $\ep\in\R_{>0}$, $n,k\in\N$ and self-adjoint contractions
$A_1,\dots,A_n\in\mathcal F$, there are $s\in\N$ and self-adjoint contractions
$B_1,\dots,B_n\in\C^{s\times s}$ such that
$$\Big|\ta(w(A_1,\dots,A_n))-\frac 1s\tr(w(B_1,\dots,B_n))\Big|<\ep
\qquad\text{for all $w\in\ax_k$.}$$
\end{enumerate}
\end{proposition}

\begin{proof}
The implication \eqref{allultra}$\Rightarrow$\eqref{oneultra} is
trivial.

For the proof of \eqref{oneultra}$\Rightarrow$\eqref{nonplusultra},
let $f\in\C\ax$ satisfy condition \eqref{algcon-matrix} from Conjecture
\ref{algcon}. Then $\ta_0(f(A_1,\dots,A_n))\ge 0$ for all
self-adjoint contractions $A_1,\dots,A_n\in\mathcal R$.
Let $\om$ be an ultrafilter on $\N$.
By \eqref{oneultra}, it suffices to show that
$\ta_{0,\om}(f(A_1,\dots,A_n))\ge 0$ for all
self-adjoint contractions $A_1,\dots,A_n\in\mathcal R^\om$.
By continuity, we may even assume that the $A_i$ are not only contractions
but there exists $\ep\in\R_{>0}$ such that $\|A_i\|\le 1-\ep$. Then
each $A_i$ has a representative $(A_i^{(j)}+B_i^{(j)})_{j\in\N}$ such that
each $A_i^{(j)}$ is a self-adjoint contraction in $\mathcal R$ and
$(B_i^{(j)})_{j\in\N}\in I_\om$. But then
\begin{align*}
\ta_{0,\om}(f(A_1,\dots,A_n))&=
\lim_{j\to\om}\ta_0(f(A_1^{(j)}+B_1^{(j)},\dots,A_n^{(j)}+B_n^{(j)}))\\
&=\lim_{j\to\om}\ta_0(f(A_1^{(j)},\dots,A_n^{(j)}))\ge 0
\end{align*}
where the second equality follows from the fact that $I_\om$ is an
ideal and $\ta_{0,\om}|_{I_{\om}}=0$.

To prove \eqref{nonplusultra}$\Rightarrow$\eqref{approximate},
let $\ep>0$ and $n,k\in\N$ be given. Consider the
finite-dimensional $\C$-vector space $\C\ax_k$ and its dual space
$\C\ax^\vee_k$. Let $C\subseteq\C\ax_k^\vee$ denote the closure of
the convex hull of the set $T\subseteq\C\ax_k^\vee$ of all the
linear forms
$$p\mapsto\frac 1s\tr(p(\bar B))\quad
(\text{$s\in\N$, $\bar B$ an $n$-tuple of self-adjoint contractions
in $\C^{s\times s}$}).$$ Now let an $n$-tuple $\bar A$ of
self-adjoint contractions in $\mathcal F$ be given and consider
$L\in\C\ax_k^\vee$ given by $L(p)=\ta(p(\bar A))$ for $p\in\C\ax_k$.

Assume $L\notin C$. By the complex Hahn-Banach separation theorem,
we then find $f\in\C\ax_k\cong\C\ax_k^{\vee\vee}$ and $c\in\R$ such
that $\re(L(f))<c<\re(L'(f))$ for all $L'\in C$. Replacing $f$ by
$f-c$, we may assume $c=0$. Then
$L'(f+f^*)=L'(f)+L'(f)^*=2\re(L'(f))>0$ for all $L'\in C$ but
$L(f+f^*)<0$, contradicting \eqref{nonplusultra}.

Therefore $L\in C$, i.e., every neighborhood of $L$ in
$\C\ax_k^\vee$ contains a convex combination of elements of $T$.
Since $\Q$ is dense in $\R$, every such neighborhood also contains
such a convex combination with rational coefficients. But building
matrices in block diagonal form, it is easy to see that the set $T$
is closed under such rational convex combinations.

To prove \eqref{approximate}$\Rightarrow$\eqref{allultra},  let
$A_1,A_2,\dots$ be a sequence of self-adjoint contractions of
$\mathcal F$ generating $\mathcal F$ as a von Neumann algebra. For
each $k\in\N$, choose self-adjoint contractions
$B_1^{(k)},\dots,B_k^{(k)}\in\mathcal R$ satisfying
$$\big|\ta(w(A_1,\dots,A_k))-
\ta_{0,\om}(w(B_1^{(k)},\dots,B_k^{(k)}))\big|<\frac
1k\qquad\text{for each $w\in\langle X_1,\dots,X_k\rangle_k$.}$$ For
each $i\in\N$, let $B_i\in\mathcal R^\om$ be the self-adjoint
contraction represented by the sequence $(B_i^{(k)})_{k\in\N}$ (with
$B_i^{(k)}:=1$ for $i>k$). Then for all $n\in\N$ and $w\in\ax$ we
have
\begin{equation}\label{limapprox}
\ta_{0,\om}(w(B_1,\dots,B_n))=
\lim_{k\to\om}\tau_0(w(B_1^{(k)},\dots,B_n^{(k)}))=\ta(w(A_1,\dots,A_n)).
\end{equation}
There is a map $\io$ that embeds the $*$-algebra
generated by the $A_i$ into $\mathcal R^\om$ by mapping $A_i$ to $B_i$
for $i\in\N$. Indeed, if $A:=\sum_w\lambda_w w(A_1,\dots,A_n)=0$
and $B:=\sum_w\lambda_w w(B_1,\dots,B_n)$, then
\eqref{limapprox} shows that $\|A\|_2=\|B\|_2$. In particular,
$\|A\|_2=0\Leftrightarrow\|B\|_2=0$ which shows that $\io$ is well-defined and
injective. By \eqref{limapprox}, it
is a trace-preserving $*$-homomorphism and therefore extends to an
embedding $\iota:\mathcal F\hookrightarrow \mathcal R^\om$.
\end{proof}

\begin{theorem}\label{nonhadwin}
The following are equivalent:
\begin{enumerate}[\rm(i)]
\item\label{non-con}
Connes' embedding conjecture {\rm \ref{concon}} holds;
\item\label{non-c}
For $\K=\C$, conditions \eqref{algcon-matrix} from Conjecture
{\rm\ref{algcon}} and the conditions from Theorem
{\rm\ref{II_1}} are equivalent for all $n\in\N$ and $f\in\C\ax$;
\item\label{non-r}
For $\K=\R$
conditions \eqref{algcon-matrix} from Conjecture {\rm\ref{algcon}} and
the conditions from Corollary {\rm\ref{contr}} are equivalent
for all $n\in\N$ and $f\in\sym\R\ax$.
\end{enumerate}
\end{theorem}

\begin{proof}
First note that condition (i) from Conjecture {\rm\ref{algcon}}
follows from the other conditions mentioned by Theorem \ref{II_1}
and Corollary \ref{contr}. Now Proposition \ref{approx} shows that
\eqref{non-con} and \eqref{non-c} are equivalent. Finally, the
equivalence of \eqref{non-c} and \eqref{non-r} follows from Proposition
\ref{cr-realtrace} together with Lemmas \ref{rc} and \ref{cr}.
\end{proof}

Combining Theorem \ref{nonhadwin} with Theorem \ref{II_1} and Corollary
\ref{contr}, we get the desired proof of Theorem \ref{mthm}.

\section{Polynomials in two variables}\label{twovarsect}

In this section, we let $n=2$ and write $(X,Y)$ instead of $(X_1,X_2)$.
Moreover, we denote by $\pi:\C\axy\to\C\cxy$ the canonical ring
epimorphism that lets the variables commute.

\begin{definition} We call a word $w\in\axy$ \emph{cyclically sorted}
if it is cyclically equivalent to $X^iY^j$ for some $i,j\in\N_0$.
A polynomial $f\in\C\axy$ is called \emph{cyclically sorted} if it
is a linear combination of cyclically sorted words.
\end{definition}

\begin{proposition}\label{csorted}
Let $f\in\C\axy$ be cyclically sorted.
Suppose $\pi(f)\ge 0$ on $[-1,1]^2$. Then
$f+\ep$ is cyclically equivalent to an element of
$M_\R$ for all $\ep\in\R_{>0}$.
\end{proposition}

\begin{proof}
For each $g\in\C\cxy$, there is exactly one linear combination $\rh(g)$
of words of the form $X^iY^j$ ($i,j\in\N_0$) such that $\pi(\rh(g))=g$.
If $p,q\in\C\axy$ are cyclically sorted and satisfy
$\pi(p)=\pi(q)$, then $p\csim q$.
The hypothesis $\pi(f)\ge 0$ on $[-1,1]^2$
implies $\pi(f)\in\R\cxy$ since the coefficients of $f$
are essentially higher partial derivatives of $f$ at the origin.
Given $\ep\in\R_{>0}$, it follows from Putinar's Theorem \ref{putinar} that
$$\pi(f)+\ep=\sum_ip_i^2+\sum_iq_i^2(1-X^2)+\sum_ir_i^2(1-Y^2)$$
for some $p_i,q_i,r_i\in\R\cxy$. This implies
$$f+\ep\csim\sum_i\rh(p_i)^*\rh(p_i)+
            \sum_i\rh(q_i)^*(1-X^2)\rh(q_i)+
            \sum_i\rh(r_i)(1-Y^2)\rh(r_i)^*\in M_\R$$
because the expressions on both sides are cyclically sorted.
\end{proof}

\begin{example}\label{first-example}
Set $$f:=(1-X^2)(1-Y^2)\in\R\axy.$$ Then $f+\ep$ is cyclically
equivalent to an element in $M_\R$ for every $\ep\in\R_{>0}$. While this
follows from Proposition \ref{csorted}, it can also be seen directly:
We may assume $\ep=\frac 1m$ for some $m\in\N$ and note that
$$f+\frac 1m\csim\Big(1-X^2+\frac 1mX^{2m}\Big)(1-Y^2)+
\frac 1m(X^mY^2X^m+(1-X^{2m})).$$ The second term of this sum lies
in $M_\R$ since
$$1-X^{2m}=\sum_{k=0}^{m-1}X^k(1-X^2)X^k,$$
and we use Remark
\ref{identity} to see that the first term is cyclically equivalent to
$$\frac 1m(1-Y^2)+\frac 1m(1-X^2)\Big(\sum_{k=0}^{m-2}(m-1-k)X^k(1-Y^2)X^k\Big)
(1-X^2)\in M_\R.$$
For $\ep=0$, $f+\ep$ is not cyclically equivalent to an element of $M_\R$.
In fact, it is an easy exercise to show that $\pi(f)\notin\pi(M_\R)$.
\end{example}

\begin{example}
The polynomial $$f:=YX^4Y+XY^4X-3XY^2X+1\in\sym\R\axy$$
is a noncommutative cyclically sorted version of the Motzkin polynomial
$\pi(f)$. The Motzkin polynomial is probably the most well-known example of
a polynomial which is nonnegative on $\R^2$ but not a sum of squares of
polynomials \cite{rez}.
By Proposition \ref{csorted}, $f+\ep$ is for each
$\ep\in\R_{>0}$ cyclically equivalent to an element of $M_\R$. This shows
in particular that $\tr(f(A,B))\ge 0$ for all $s\in\N$ and all
self-adjoint contractions $A,B\in\C^{s\times s}$. Since $\pi(f)\ge 0$ on
(any square in) $\R^2$, we can use the same reasoning together with a scaling
argument to see that
$\tr(f(A,B))\ge 0$ for all $s\in\N$ and all
self-adjoint matrices $A,B\in\C^{s\times s}$, a fact for which we do not
know a direct proof. However, a direct proof that $f+\ep$ is for all
$\ep\in\R_{>0}$ cyclically equivalent to an element of $M_\R$ can be obtained
as in the previous example since
$$f\csim Y(1-X^2)^2Y+X(1-Y^2)^2X+(1-X^2)(1-Y^2).$$
Note that $f(A,B)$ is not positive semidefinite for all
self-adjoint contractions $A,B\in\R^{2\times 2}$, since for
$$
A:=\frac 12
\begin{pmatrix}
1&1\\
1&1
\end{pmatrix}
\qquad\text{and}\qquad
B:=\begin{pmatrix}
-1&0\\
0&1
\end{pmatrix},
$$
$$f(A,B)=\frac 12 \begin{pmatrix}
\hfill 1&-3\\
-3&\hfill 1
\end{pmatrix}
$$
is clearly not positive semidefinite.
\end{example}

\section{Bounds}\label{bousect}
In this section, we use valuation theory \cite{pc}, basic first
order logic and model theory of real closed fields \cite{pre} to
derive certain bounds for Conjecture \ref{algcon}. For the moment,
let \eqref{algcon-matrix} and \eqref{algcon-module} refer to the
respective conditions for $\K=\R$ in Conjecture \ref{algcon}. As we
have seen in Theorem \ref{mthm}, Connes conjecture is equivalent to
the implication
\eqref{algcon-matrix}$\Rightarrow$\eqref{algcon-module} for
$f\in\sym\R\ax$. Here we show that this implication must actually
hold in a stronger form if it holds at all. Suppose that Connes'
conjecture holds and we are given $f\in\sym\R\ax$ and
$\ep\in\R_{>0}$. Then there are two bounds. First, there is a bound
on the size of the matrices on which the nonnegativity condition in
\eqref{algcon-matrix} has to be tested. Second, there is a bound on
the degree complexity of the representation of $f+\ep$ (for this
particular $\ep$) in \eqref{algcon-module}. These bounds depend only
on $\ep$, the number of variables, the degree of $f$ and the size of the
coefficients of $f$ (rather than on $f$ itself).
Moreover, the bounds are computable from this data (in the sense of
recursion theory). Unfortunately, the rather nonconstructive methods
yielding these bounds do not allow for further specification of the
kind of dependence. We will first prove a certain technical version
of Corollary \ref{contr} which is valid not only over $\R$ but over
any real closed field (see Proposition \ref{rcf}).
\smallskip

Let us recall some facts from the theory of ordered fields.
Suppose $R$ is a real closed field. Let $\le$ denote the ordering of $R$ and
$$\O:=\{a\in R\mid|a|\le N\text{\ for some $N\in\N$}\}$$
the convex hull of $\Z$ in $R$. This is a valuation ring with
(unique) maximal ideal $\m$ given by
$$\m=\{a\in R\mid N|a|\le 1\text{\ for all $N\in\N$}\}.$$
The residue field $\O/\m$ is again a real closed field
(cf.~\cite[8.6]{pre} or \cite[II \S4 Lemma 17]{pc}), but this
time archimedean and thus embeds uniquely into $\R$
\cite[II \S 3 Satz 3]{pc}. We therefore
always assume $\O/\m\subseteq\R$. Moreover, we find at least one embedding
$\rh:\O/\m\hookrightarrow\O\subseteq R$ such that
$\overline{\rh(x)}=x$ for all $x\in\O/\m$ \cite[III \S 2 Satz 6]{pc}.
We extend the canonical
homomorphism $\O\to\O/\m\subseteq\R$ to a ring homomorphism
$$\O\ax\to\R\ax,\ f\mapsto\overline f$$
sending $X_i$ to $X_i$. Similarly, $\rh$ can be extended to polynomials.

The quadratic module
$M_R\subseteq\sym R\ax$ generated by $1-X_1^2,\dots,1-X_n^2$
consists exactly of the sums of elements of the form
\begin{equation}\label{ofthisform}
g^*g\qquad\text{and}\qquad g^*(1-X_i^2)g\qquad
(1\le i\le n,\ g\in R\ax).
\end{equation}
Now consider only elements of this form
of degree at most $2k$ ($k\in\N$) and call the set of all sums of such
elements $M_{R,k}$. Then $M_{R,k}$ is a convex cone in the $R$-vector space
$\sym R\ax_{2k}$ which is
(perhaps strictly) contained in $M_R\cap R\ax_{2k}$.
Clearly, $M_R=\bigcup_{k\in\N}M_{R,k}$.

Since we will no longer be concerned with complex matrices but with matrices
over real closed fields, it seems more appropriate to speak of \emph{symmetric}
matrices rather than self-adjoint ones.

\begin{proposition}\label{closed}
Suppose $k\in\N$. Let $U$ denote the subspace of $\sym\R\ax_{2k}$
of those elements which are cyclically equivalent to $0$. Then $M_{\R,k}+U$
is closed in $\sym\R\ax_{2k}$.
\end{proposition}

\begin{proof}
Let $\pi:\sym\R\ax_{2k}\to(\sym\R\ax_{2k})/U=:V$ be the canonical
projection. Then $M_{\R,k}+U=\pi^{-1}(\pi(M_{\R,k}))$. Hence, it
suffices to show that the convex cone $\pi(M_{\R,k})$ is closed in
$V$. By Carath\'eodory's theorem (see e.g. \cite[p.~40, Exercise
1.8]{hol}), each element of $\pi(M_{\R,k})$ can be written as the
image of a sum of at most $m$ terms of the form \eqref{ofthisform}
where $m:=\dim V$. Setting $p_0:=1$ and $p_i:=1-X_i^2$ for
$i\in\{1,\dots,n\}$, we see that $\pi(M_{\R,k})$ is the image of the
map
$$\Ph:\begin{cases}
      \R\ax_k^m\times\R\ax_{k-1}^m\times\dots\times\R\ax_{k-1}^m\to V\\
      \ (g_{01},\dots,g_{0m},\quad\dots\quad,g_{n1},\dots,g_{nm})
      \mapsto\pi\left(\sum_{i=0}^n\sum_{j=1}^mg_{ij}^*p_ig_{ij}\right).
      \end{cases}
$$
We claim that $\Ph^{-1}(0)=\{0\}$. To show this, suppose
\begin{equation}\label{s0}
h:=\sum_{i=0}^n\sum_{j=1}^mg_{ij}^*p_ig_{ij}\csim 0.
\end{equation}
Let $s\in\N$ and $A_1,\dots,A_n\in\R^{s\times s}$ be symmetric with
$\|A_i\|<1$. Then $1-A_i^2$ is a positive definite
and can be written as $1-A_i^2=B_i^2$ for some
symmetric invertible $B_i\in\R^{s\times s}$. It is convenient to let
$B_0$ denote the identity matrix in $\R^{s\times s}$. Denoting by $e_t$ the
$t$-th unit vector of $\R^s$, it follows from \eqref{s0} that
\begin{align*}
\sum_{t=1}^s\sum_{i=0}^n\sum_{j=1}^m\langle
B_ig_{ij}(\bar A)e_t,
B_ig_{ij}(\bar A)e_t\rangle=
\tr(h(A_1,\dots,A_n))=0.
\end{align*}
Consequently, we get $B_ig_{ij}(\bar A)e_t=0$ and hence
$g_{ij}(\bar A)e_t=0$ for all $i,j,t$. This shows that
$g_{ij}(A_1,\dots,A_n)=0$ for all symmetric $A_i\in\R^{s\times s}$ with
$\|A_i\|<1$. By continuity, the same holds for all symmetric
contractions
$A_i\in\R^{s\times s}$. Hence \cite[Proposition 2.3]{ks} implies
that $g_{ij}=0$. This shows that $\Ph^{-1}(0)=\{0\}$. Together with the fact
that $\Ph$ is homogeneous, \cite[Lemma 2.7]{ps} shows that
$\Ph$ is a proper and therefore a closed map. In particular, its image
$\pi(M_{\R,k})$ is closed in $V$.
\end{proof}

In the following lemma, we will apply \emph{Tarski's transfer principle}, i.e.,
the
fact that exactly the same first order sentences with symbols $0$, $1$, $+$,
$\cdot$, $\le$ hold in each real closed field \cite[5.3]{pre}.

\begin{lemma}\label{sepp}
Let $k\in\N$ and $U$ be the subspace of $\sym R\ax_{2k}$
of those elements which are cyclically equivalent to $0$.
Suppose that $f\in\sym R\ax_{2k}\setminus(M_{R,k}+U)$. Then there is a linear
map $L:\sym R\ax_{2k}\to R$ such that $L(M_{R,k})\subseteq R_{\ge 0}$,
$L|_U=0$, $L(1)=1$ and $L(f)<0$.
\end{lemma}

\begin{proof}
We first prove this for $R=\R$. Consider the convex cone
$M_{\R,k}+U$ in $\sym\R\ax_{2k}$ which is closed by Proposition
\ref{closed}. Separating this cone from the cone spanned by a little
ball around $f$ (use e.g. \cite[p.~15, \S 4B Corollary]{hol}), we
find a linear map $L_0:\sym\R\ax_{2k}\to\R$ such that
$L_0(M_{\R,k}+U)\subseteq\R_{\ge 0}$ and $L_0(f)<0$. Since $1\in
M_{\R,k}$, we have $L_0(1)\ge 0$. If $L_0(1)>0$, then
$L:=\frac{L_0}{L_0(1)}$ has the desired properties. If $L_0(1)=0$,
then we set $L:=L_1+\la L_0$ where
$$L_1:\sym\R\ax_{2k}\to\R,\quad g\mapsto g(0)$$ and
$\la\in\R_{>0}$ is sufficiently large to ensure that $L(f)<0$.
This proves the statement for $R=\R$.

The general case follows by Tarski's transfer principle once we know that the
statement can \emph{for fixed} $k,n\in\N$ be expressed in the first order
language with symbols $0$, $1$, $+$, $\cdot$, $\le$. But this is indeed
possible: To model $f\in\sym R\ax_{2k}$, use universal
quantifiers for the finitely many coefficients that a polynomial of degree
$2k$ in $n$ variables can have. The condition $f\notin M_{R,k}+U$ can also
be written down in this language by using Carathéodory's theorem as in the
proof of Proposition \ref{closed}. The existence of the linear map $L$ can be
expressed by existential
quantifiers for the values of $L$ on a basis of $\sym R\ax_{2k}$.
\end{proof}

By Lemma \ref{algint1} and \eqref{non-useful}, we find for every
word $w\in\ax$ an $N_w\in\N$ such that
$N_w\pm(w+w^*)\in M_\Q$. Moreover, we find for each $k\in\N$
some $d_k\ge k$ such that
\begin{equation}\label{dk}
2N_w\pm(w+w^*)\in M_{\Q,d_k}\subseteq M_{R,d_k}\qquad
\text{for all $w\in\ax_{2k}$.}
\end{equation}

\begin{lemma}\label{duality}
Suppose $k\in\N$ and $f\in\sym R\ax_{2k}$ is not cyclically
equivalent to an element of $M_{R,k}$. Then there is a linear map
$L:R\ax_{2k}\to R$ such that $L(f)<0$,
\begin{enumerate}[\rm (a)]
\item $L(pq)=L(qp)$ for all $p,q\in R\ax$ such that $pq\in R\ax_{2k}$;
\label{du-trace}
\item $L(M_{R,k})\subseteq R_{\ge 0}$;
\label{du-m}
\item $|L(w)|\le N_w$ for all $w\in\ax_{2k}$;
\label{du-bound}
\item $L(1)=1$;
\label{du-1}
\item $L(p^*)=L(p)$ for all $p\in R\ax_{2k}$.
\label{du-*}
\end{enumerate}
\end{lemma}

\begin{proof}
Set $d:=d_k\ge k$. By Lemma \ref{sepp}, we find a linear map
$L_0:\sym R\ax_{2d}\to R$ such that $L_0(M_{R,d})\subseteq R_{\ge 0}$,
$L_0|_U=0$, $L_0(1)=1$ and $L_0(f)<0$ where $U\subseteq\sym R\ax_{2d}$ is the
subspace of
polynomials that are cyclically equivalent to $0$. The linear map
$$L:R\ax_{2k}\to R,\quad p\mapsto L_0\left(\frac{p+p^*}2\right)$$
extends the restriction of $L_0$ to $\sym R\ax_{2k}$ which shows
\eqref{du-m},\eqref{du-1} and $L(f)<0$. Property \eqref{du-*} is
clear from the definition of $L$. By \eqref{dk}, we have
$$2(N_w\pm L(w))=2N_wL(1)\pm(L(w)+L(w^*))=
L_0(2N_w\pm(w+w^*))\ge 0$$
which yields \eqref{du-bound}. To show \eqref{du-trace}, suppose
$p,q\in R\ax$ are such that $pq\in R\ax_{2k}$. Then $pq\csim qp$ and
$(pq)^*\csim(qp)^*$ imply that $pq+(pq)^*\csim qp+(qp)^*$. This shows
$pq+(pq)^*-(qp+(qp)^*)\in U$ whence
$2L(pq)=L_0(pq+(pq)^*)=L_0(qp+(qp)^*)=2L(qp)$.
\end{proof}

\begin{lemma}\label{pushdown}
Suppose $k\in\N$ and $f\in\sym\O\ax_{2k}$ is not cyclically
equivalent to an element of $M_{R,k}$. Then there is a linear map
$L:(\O/\m)\ax_{2k}\to\O/\m$ that satisfies $L(\overline f)\le 0$
and conditions \eqref{du-trace}--\eqref{du-*} from Lemma {\rm\ref{duality}}
{\rm (}with $R$ replaced by $\O/\m${\rm )}.
\end{lemma}

\begin{proof}
Let $L_0$ be one of the linear maps whose existence has been shown in
the previous lemma. Property \eqref{du-bound}
(with $L$ replaced by $L_0$) implies that $L_0(\O\ax)\subseteq\O$. We can thus
define the map
$$L:(\O/\m)\ax_{2k}\to\O/\m,\quad p\mapsto\overline{L_0(\rh(p))}.$$
Using that $\overline{\rh(\la)}=\la$ for all $\la\in\O/\m$, we see that
$L$ is $\O/\m$-linear. We know that
$\rh(\overline f)-f$ has all its coefficients in $\m$. Because of property
\eqref{du-bound}, this shows that $L_0(\rh(\overline f)-f)\in\m$
whence
$$L(\overline f)=\overline{L_0(\rh(\overline f))}
=\overline{L_0(f)}+\overline{L_0(\rh(\overline
f)-f)}=\overline{L_0(f)}\le 0.$$ Moreover, it is easy to see that
$L$ inherits properties \eqref{du-trace}--\eqref{du-*} from $L_0$.
\end{proof}

\begin{lemma}\label{pullup}
Suppose $k\in\N$ and $f\in\sym\O\ax_{2k}$ is not cyclically
equivalent to an element of $M_{R,k}$. Then there is a linear map
$L:\R\ax_{2k}\to\R$ that satisfies $L(\overline f)\le 0$
and conditions \eqref{du-trace}--\eqref{du-*} from Lemma {\rm\ref{duality}}
{\rm (}with $R$ replaced by $\R${\rm )}.
\end{lemma}

\begin{proof}
Let $L_0$ be one of the linear maps whose existence has been shown in
the previous lemma.
Let $x_w$ and $y_w$ be variables in the
formal language of first order logic where $w$ ranges over all $w\in\ax_{2k}$.
Build up a formula $\Ph$ with free variables $x_w$ and
$y_w$ in the first order language with symbols
$0$, $1$, $+$, $\cdot$, $\le$ expressing that (over the real closed field
$R$ where the formula is interpreted) $L(\sum_w y_ww)\le 0$ and conditions
\eqref{du-trace}--\eqref{du-*} from Lemma \ref{duality}
hold for the linear map $L:R\ax_{2k}\to R$ given by $L(w)=x_w$.
Compare the second part of the proof of Lemma \ref{sepp}
for some details on how this can be done. By Lemma \ref{pushdown}, $\Ph$ holds
in the real closed field $\O/\m$ when $x_w$ is interpreted as $L_0(w)$ and
$y_w$ is interpreted as the coefficient of $w$ in $\overline f$.
Define another formula
$\Ps$ with free variables $y_w$ which arises from $\Ph$ by quantifying
all $x_w$ existentially. Then $\Ps$ holds in $\O/\m$ when the $y_w$ are
interpreted as the coefficients of $f$. By the substructure completeness of
the theory of real closed fields \cite[5.1,4.7]{pre},
$\Ps$ holds also in the real closed extension field $\R$ of $\O/\m$
under the same interpretation of the $y_w$.
\end{proof}

\begin{lemma}\label{perturbation}
Suppose $f\in\sym R\ax_{2k}$ has all its coefficients in $\m$. Then for each
$\ep\in R_{>0}\setminus\m$, we have $f+\ep\in M_{R,d_k}$.
\end{lemma}

\begin{proof}
Without loss of generality, we may assume that
$f=a(w+w^*)$ with $a\in\m$ and $w\in\ax_{2k}$. Then
\begin{align*}
f+\ep&=a(w+w^*)+|a|N_w+(\ep-|a|N_w)\\
&=|a|(N_w+\sign(a)(w+w^*))+(\ep-|a|N_w)\in M_{R,d_k}
\end{align*}
since $\ep-|a|N_w\ge 0$ and $N_w\pm(w+w^*)\in M_{\Q,d_k}\subseteq M_{R,d_k}$
by \eqref{dk}.
\end{proof}

\begin{proposition}\label{rcf}
Suppose $f\in\sym\O\ax$ and $\ph(\overline f)\ge 0$
for all tracial contraction states $\ph$ on $\R\ax$. Then for all
$\ep\in R_{>0}\setminus\m$, $f+\ep$ is cyclically equivalent to an element of
$M_R$.
\end{proposition}

\begin{proof}
We show the contraposition, i.e., we assume that we have
$N\in\N$ such that $f+\frac 1N$ is \emph{not}
cyclically equivalent to an element of $M_R$ and find a tracial
contraction state $\ph$ on $\R\ax$ such that $\ph(\overline f)<0$.
Let \eqref{du-trace}--\eqref{du-*} refer to the conditions from Lemma
\ref{duality} with $R$ replaced by $\R$. Lemma \ref{pullup} provides us
for each $k\in\N$ such that $2k\ge\deg f$ with a linear map
$L_k:\R\ax_{2k}\to\R$ satisfying $L_k(\overline f+\frac 1N)\le 0$ and
\eqref{du-trace}--\eqref{du-*}. To each $L_k$, we associate a point $P_k$ in
the product space $S:=\prod_{w\in\ax}[-N_w,N_w]$
by setting $P_k(w):=L_k(w)$ if $w\in\ax_{2k}$ and
$P_k(w):=0$ if $w\in\ax\setminus\ax_{2k}$. Since $S$ is compact
by Tychonoff's theorem, the sequence $(P_k)_k$ has a
subsequence converging to some $P\in S$. Define the linear map
$\ph:\R\ax\to\R$ by $\ph(w):=P(w)$ for all $w\in\ax$. Using \eqref{du-m},
\eqref{du-1} together with $M_R=\bigcup_{k\in\N}M_{R,k}$,
\eqref{du-trace}, \eqref{du-*}
and Lemma \ref{tcs-module}, it is easy to see that $\ph$ is a tracial
contraction state such that $\ph(\overline f+\frac 1N)\le 0$ and
therefore $\ph(\overline f)\le-\frac 1N<0$.
\end{proof}

\begin{theorem}\label{boundthm}
Suppose that Connes' embedding conjecture {\rm \ref{concon}} holds. Then there
is a computable function $N:\N\to\N$ such that for all $t\in\N$ the following
is true: Whenever $n\in\N$ with $n\le t$,
$f\in\sym\R\langle X_1,\dots,X_n\rangle$ is of degree $\le t$, has
absolute value of its coefficients bounded by $t$ and
satisfies $\tr(f(A_1,\dots,A_n))\ge 0$ for all symmetric contractions
$A_i\in\R^{N(t)\times N(t)}$, then $f+\frac 1t$ is cyclically equivalent to an
element of $M_{\R,N(t)}$.
\end{theorem}

\begin{proof} For technical reasons, it is convenient to replace the condition
$n\le t$ in the statement by the condition $n=t$. This does not affect the
generality of the theorem
since $$M_{\R,N}^{(t)}\cap\R\ax=M_{\R,N}^{(n)}\qquad\text{for $n\le t$}.$$
Most facts about finite-dimensional real Euclidean vector spaces
carry over from $\R$ to any real closed field by Tarski's transfer principle.
We will
therefore use concepts like symmetric contractions
over the real closed fields $R$ and $\O/\m$. For a matrix $A\in\O^{s\times s}$,
we can apply the map $\O\to\O/\m$ entrywise and get a matrix
$\overline A\in(\O/\m)^{ s\times s}$. For every
symmetric contraction $A\in(\O/\m)^{s\times s}$, there
is a symmetric contraction $B\in R^{s\times s}$ with all its entries in $\O$
such that $\overline B=A$.

{\sc Claim 1.} For fixed $t\in\N$, the following infinitely many conditions
(a), (b), (c), (d$_s$) and (e$_s$) ($s\in\N$) cannot be satisfied
simultaneously.
\begin{enumerate}[(ei)]
\item[(a)] $R$ is a real closed field;
\item[(b)] $f\in\sym R\langle X_1,\dots,X_t\rangle$ is of degree at most $t$;
\item[(c)] The absolute value of the coefficients of $f$ is bounded by $t$;
\item[(d$_s$)] $\tr(f(A_1,\dots,A_t))\ge 0$ for all symmetric contractions
$A_i\in R^{s\times s}$;
\item[(e$_s$)] $f+\frac 1t$ is not cyclically equivalent to an element of
$M_{R,s}^{(t)}$.
\end{enumerate}
\emph{Proof of Claim {\rm 1}.} Assuming these conditions, we obtain the
following.
\begin{enumerate}[(ei)]
\item[(c$'$)] $f\in\O\langle X_1,\dots,X_t\rangle$;
\item[(b$'$)] $\overline f\in\sym\R\langle X_1,\dots,X_t\rangle$;
\item[(d$_s'$)] $\tr(\overline f(A_1,\dots,A_t))\ge 0$ for all symmetric
contractions $A_i\in\R^{s\times s}$;
\item[(e$'$)] $\overline f+\frac 1{2t}$ is not cyclically equivalent to an
element of $M_{\R}^{(t)}$.
\end{enumerate}
Of course, (c$'$) follows from (c) by the definition of
$\O$. Because of (c$'$), we can consider $\overline f$ and from (b)
it is clear that (b$'$) holds. It is easy to see that (d$_s$) implies (d$_s'$)
for all
symmetric contractions $A_i\in(\O/\m)^{s\times s}$. With Tarski's transfer
principle
and the fact that $\O/\m$ and $\R$ are real closed, it is easy to extend this
from $\O/\m$ to $\R$ (cf. Lemma \ref{pullup}). Now assume that (e') does not
hold, i.e.,
$\overline f+\frac 1{2t}$ is cyclically equivalent to an element of
$M_{\R,s}^{(t)}$ for some $s\in\N$. By Tarski's principle (use again
Carath\'eodory's theorem to express this in first order logic), we get
\begin{equation}\label{mrseq}
\rh(\overline f)+\frac 1{2t}\text{\ is cyclically equivalent to an element
of\ }M_{R,s}^{(t)}\subseteq M_{R,d_s}^{(t)}.
\end{equation}
From the fact that $f-\rh(\overline f)$ has all its coefficients in $\m$
and Lemma \ref{perturbation}, it follows that
$f-\rh(\overline f)+\frac 1{2t}\in M_{R,d_s}^{(k)}$. Combining this with
\eqref{mrseq} yields that $f+\frac 1t$ is cyclically equivalent to an
element of $M_{R,d_s}^{(k)}$ which contradicts
$(e_{d_s})$. Finally, use Proposition \ref{rcf} to see
that (b$'$), (d$_s'$) ($s\in\N$) and
(e$'$)
cannot be satisfied simultaneously if the algebraic version \ref{algcon} of
Connes' conjecture holds. But this algebraic version is equivalent to
Connes' conjecture by Theorem \ref{mthm}. This proves Claim 1.
\smallskip

As we have just seen, a lot of specifications (like the degree in (b), the
concrete bound for the absolute value of the coefficients in (c), etc.) are
not needed for Claim 1 but they ensure that the next claim holds.

{\sc Claim 2.} For fixed $t\in\N$, the above conditions (a), (b), (c), (d$_s$)
and (e$_s$)
($s\in\N$) can be expressed in the language of first order logic with
symbols $0$, $1$, $+$, $\cdot$, $\le$ and new constants for the
finitely many coefficients that a polynomial
$f\in R\langle X_1,\dots,X_t\rangle$ of degree at most $t$ can have.
Moreover, there is a \emph{decidable} (i.e., recursive) set of formulas in this
language corresponding to (a), (b), (c), (d$_s$) and (e$_s$).

\emph{Proof of Claim {\rm 2}.}
Concerning (a), write down the axioms for real closed fields.
For (b), we have introduced the new constants.
The natural number $t$ in (c) can be
written as $1+\dots+1$. There are several good ways to express (d$_s$) by a
formula for each fixed $s$. Finally, use Carath\'eodory's
theorem once more to translate (e$_s$) into such a formula for each fixed $s$.

{\sc The algorithm.} We describe a procedure how to calculate the function
$N$ that we are looking for. The program takes $t\in\N$ and
yields a suitable $N(t)$. Let the program generate
successively all words of length $1,2,3,\dots$ over the finite alphabet of the
language from Claim 1. Every time a word has been generated, let the program
check whether this is by chance a formal proof of $0=1$ in the first order
predicate calculus that uses only axioms from the set of formulas from Claim 1
(this can be checked since this set is decidable by Claim 2).
When the program encounters such a formal proof, let it terminate after
outputting the smallest number $N(t)$ such that the found formal proof uses
as axioms only (a), (b), (c), (d$_s$) and (e$_s$) for $s\le N(t)$.

\emph{Proof of termination.} Since the set of allowed axioms is inconsistent
by Claim 1, $0=1$ is a logical consequence of it. By Gödel's completeness
theorem, the algorithm will thus eventually terminate.

\emph{Proof of correctness.} The number $N(t)$ has the desired properties
because $\R$ is real closed and conditions (a), (b), (c), (d$_s$), (e$_s$)
for $s:=N(t)$ must be inconsistent (observe that (d$_{k+1}$) implies (d$_k$)
and (e$_{k+1}$) implies (e$_k$) for all $k\in\N$).
\end{proof}

Note that the information that the bound $N(t)$ is computable from
$t$ means that it can in a certain sense not grow ``too'' fast when
$t\to\infty$. By a diagonal argument, it is indeed easy to
see that there are functions $\N\to\N$ growing faster than any
computable function. On the other hand, the described algorithm
computing $N(t)$ from $t$ has a tremendous complexity and is
therefore purely theoretical. If one is not interested in the
information that $N$ is computable, one can replace Gödel's
completeness theorem by the compactness theorem from first order
logic.

\subsection*{Acknowledgments.} The second author wants to thank Prof. Dr.
J. Cimpri\v c for the invitation to Ljubljana where part of this work was done.
Both authors would like to thank Prof. Dr. M. Putinar for introducing them to
R\u adulescu's work \cite{r2} on Connes' conjecture.

\end{document}